%% file: DS.tex
\documentclass[12pt]{article}
\usepackage{amsmath}
\usepackage{amssymb}
\def\eps{\varepsilon}
\font\tencmmib=cmmib10 \skewchar\tencmmib '60
\newfam\cmmibfam
\textfont\cmmibfam=\tencmmib

\def\bbox{\quad\hbox{\vrule \vbox{\hrule \vskip2pt \hbox{\hskip2pt
\vbox{\hsize=1pt}\hskip2pt} \vskip2pt\hrule}\vrule}}
\def\lessim{\ \lower4pt\hbox{$
\buildrel{\displaystyle <}\over\sim$}\ }
\def\gessim{\ \lower4pt\hbox{$\buildrel{\displaystyle >}
\over\sim$}\ }
\def\P{{\cal P}}

\def\DD{{\cal D}}

\def\M{{\cal M}}

\def\eps{{\varepsilon}}

\def\id{{\mbox{\rm id}}}

\def\qed{\hfill\break\rightline{$\bbox$}}
\newcommand{\e}{\mathbb{E}}
\newcommand{\p}{\mathbb{P}}
\newcommand{\q}{\mathbb{Q}}
\newcommand{\Reals}{\mathbb{R}}

\newcommand{\F}{{\cal F}}

\newtheorem{proposition}{Proposition}
\newtheorem{lemma}{Lemma}

\makeatletter
\@addtoreset{equation}{section}

\makeatother

\input invlat.tex

\begin{document}

\title{
On the Dovbysh-Sudakov representation result.}

\author{ 
Dmitry Panchenko\thanks{
Department of Mathematics, Texas A\&M University, 
Mailstop 3386, College Station, TX, 77843,
email: panchenk@math.tamu.edu. Partially supported by NSF grant.}\\
{\it Texas A\&M University}
}
\date{}

\maketitle

\begin{abstract}
We present a detailed proof of the Dovbysh-Sudakov representation for symmetric positive 
definite weakly exchangeable infinite random arrays, called Gram-de Finetti matrices,
which is based on the representation result of Aldous and Hoover for arbitrary 
(not necessarily positive definite) symmetric weakly exchangeable arrays.
\end{abstract}
\vspace{0.5cm}

Key words: exchangeability, spin glasses.

Mathematics Subject Classification: 60G09, 82B44.

\section{Introduction.}

We consider an infinite random matrix $R= (R_{l,l'})_{l,l'\geq 1}$ which is symmetric, 
nonnegative definite in a sense that $(R_{l,l'})_{1\leq l,l'\leq n}$ is nonnegative definite 
for any $n\geq 1,$ and weakly exchangeable, which means that for any $n\geq 1$ and for any
permutation $\rho$ of $\{1,\ldots, n\}$ the matrix $(R_{\rho(l),\rho(l')})_{1\leq l,l' \leq n}$
has the same distribution as $(R_{l,l'})_{1\leq l,l' \leq n}.$ Following \cite{DS},
we will call a matrix with the above properties a Gram-de Finetti matrix.
Since all its properties - symmetric, positive definite and weakly exchangeable - 
are expressed in terms of its finite dimensional distributions, we can think of $R$ 
as a random element in the product space $M = \prod_{1\leq l,l'} \Reals$ 
with the pointwise convergence topology and the Borel  $\sigma$-algebra $\M.$ 
Let $\P$ denote the set of all probability measures on $\M.$ 
Suppose that $\p\in \P$ is such that for all $A\in \M,$
\begin{equation}
\p(A) = \int_\Omega \q(u, A) \,d\!\Pr(u) 
\label{mixture}
\end{equation}
where $\q: \Omega\times \M \to [0,1]$ is a probability kernel from
some probability space $(\Omega, \F, \Pr)$ to $M$ such that  
(a) $\q(u,\cdot)\in \P$ for all $u\in\Omega$ and
(b) $\q(\cdot, A)$ is measurable on $\F$ for all $A\in \M.$
In this case we will say that $\p$ is a mixture of laws $\q(u,\cdot)$.
We will say that a law $\q \in \P$ of a Gram-de  Finetti matrix is generated 
by an i.i.d. sample if there exists a probability measure $\eta$ on $\ell^2\times\Reals^+$ 
such that $\q$ is the law of
\begin{equation}
\bigl( h_l\cdot h_{l'} + a_l\,\delta_{l,l'}\bigr)_{l,l'\geq 1}
\label{repr}
\end{equation}
where $(h_l, a_l)$ is an i.i.d. sequence from $\eta$ and $h\cdot h'$ denotes the scalar product  
on $\ell^2$. For simplicity, we will often say that a matrix (rather than its law on $\M$) is generated
by an i.i.d. sample from measure $\eta$. The result of L.N. Dovbysh and V.N. Sudakov 
in \cite{DS} states the following.
\begin{proposition}\label{DS}
A law $\p\in\P$ of any Gram-de Finetti matrix is a mixture (\ref{mixture}) of laws in 
$\P$ such that for all $u\in\Omega$, $\q(u,\cdot)$ is generated by an i.i.d. sample.
\end{proposition}
Proposition \ref{DS} has recently found important applications in spin glasses; for example, 
it played a significant role in the proof of the main results in \cite{AA} and \cite{GG},
where a problem of ultrametricity of an infinite matrix $(R_{l,l'})_{l,l'\geq 1}$ was considered 
under various hypotheses on its distribution. For this reason, it seems worthwhile to have an 
accessible proof of this result which was, in fact, the main motivation for writing this paper.
Currently, there are two known proofs of Proposition \ref{DS}. The proof in the original paper
 \cite{DS} contains all the main ideas that will appear, maybe in a somewhat different form, 
 in the present paper but the proof is too condensed and does not provide enough details necessary 
 to penetrate these ideas. Another available proof in \cite{Hestir} is much more detailed but, unfortunately, 
it is applicable not to all Gram-de Finetti matrices even though it works in certain cases. 

In the present paper we will give a detailed proof of Proposition \ref{DS} which starts with 
exactly the same idea as \cite{Hestir}. Namely, we will deduce Proposition \ref{DS} from the 
representation result for arbitrary weakly exchangeable arrays that are not necessarily positive definite, 
due to D. Aldous (\cite{Aldous}, \cite{Aldous2}) and D.N. Hoover (\cite{Hoover}, \cite{Hoover2}),
which states that for any weakly exchangeable matrix there exist two measurable functions 
$f:[0,1]^4\to\Reals$ and $g:[0,1]^2\to\Reals$ such that the distribution of the matrix coincides with 
the distribution of 
\begin{equation}
R_{l,l} = g(u,u_l)
\,\mbox{ and }\,
R_{l,l'} = f(u,u_l,u_{l'},u_{l,l'})
\,\mbox{ for }\,
l\not = l',
\label{fg}
\end{equation}
where random variables $u,(u_l), (u_{l,l'})$ are i.i.d. uniform on $[0,1]$ and function $f$ is symmetric 
in the middle two coordinates, $u_l$ and $u_{l'}$. It is customary to define the diagonal elements as a 
function of three variables $R_{l,l} = g(u,u_l,v_l)$ where $(v_l)$ is another i.i.d. sequence with uniform
distribution on $[0,1]$; however, one can always express a pair $(u_l,v_l)$ as a function of one uniform random variable $u_l'$ in order to obtain the representation (\ref{fg}). We will consider a weakly
exchangeable matrix defined by (\ref{fg}) and, under an additional assumption that it is positive 
definite with probability one, we will prove that its distribution is a mixture of distributions generated 
by an i.i.d. sample in the sense of (\ref{repr}). First, in Section \ref{SecBounded} we will consider 
a uniformly bounded case, $|f|,|g|\leq 1,$ and then in Section \ref{SecUnbounded} we will show 
how the unbounded case follows by a truncation argument introduced in \cite{DS}. In the general
case of Section \ref{SecUnbounded} we do not require any integrability conditions on $g$ rather 
than $g<+\infty.$ Finally, to a reader interested in the proof of (\ref{fg}) we recommend a comprehensive
recent survey \cite{Austin} of the representation results for exchangeable arrays.

\textbf{Acknowledgment.} The author would like to thank Gilles Pisier and Joel Zinn for several helpful conversations.

\section{Bounded case.}\label{SecBounded}

We will start with the case when the matrix elements $|R_{l,l'}|\leq 1$ for all $l,l'\geq 1$
with probability one,  so we can assume that both functions $|f|, |g|\leq 1.$ Of course, 
the representation of the law of $R$  as the mixture (\ref{mixture}) will be simply the 
disintegration of the law of (\ref{fg}) on $\M$ over the first coordinate $u.$ The main
problem is now to show that for a fixed $u$ the (law of) matrix $R$ can be represented as (\ref{repr}).
In other words, if we make the dependence of $f$ and $g$ on $u$ implicit, then assuming that a weakly 
exchangeable matrix given by
\begin{equation}
R_{l,l} = g(u_l)
\,\mbox{ and }\,
R_{l,l'} = f(u_l,u_{l'},u_{l,l'})
\,\mbox{ for }\,
l\not = l'
\label{fg2}
\end{equation}
is positive definite with probability one, we need to show that its law can be represented as (\ref{repr}).
Our first step is to show that $f$ does not depend on the last coordinate, which is exactly the same 
as Lemma 3 in \cite{Hestir}.
\begin{lemma}\label{LemNo}
If $R$ in (\ref{fg2}) is positive definite with probability one then for
$$
\bar{f}(x,y) = \int_0^1\! f(x,y,u)\,du
$$
we have $f(u_1,u_2,u_{1,2}) = \bar{f} (u_1,u_2)$ a.s.
\end{lemma}
\textbf{Proof.} We will give a sketch of the proof for completeness. Since $(R_{l,l'})$ is positive 
definite, for any sequence of bounded measurable functions $(h_l)$ on $[0,1],$
\begin{equation}
\frac{1}{n}\sum_{l,l'=1}^n \e' R_{l,l'} h_l(u_l) h_{l'}(u_{l'})\geq 0
\label{posit}
\end{equation}
almost surely, where $\e'$ denotes the expectation in $(u_l)$. Let us take $n=4m$ 
and given two measurable sets $A_1,A_2\subset [0,1]$, let $h_l(x)$ be equal to
$$
\mbox{
$I(x\in A_1)$ for $1\leq l\leq m,$
$\,-I(x\in A_1)$ for $m+1\leq l\leq 2m,$
}
$$
$$
\mbox{
$I(x\in A_2)$ for $2m+1\leq l\leq 3m,$
$\,-I(x\in A_2)$ for $3m+1\leq l\leq 4m.$
}
$$
With this choice of $(h_l)$, the sum over the diagonal terms $l=l'$ in (\ref{posit}) is a constant,
$$
\frac{1}{2}\Bigl(\int_{A_1} \!g(x)\,dx +\int_{A_2} \!g(x)\,dx\Bigr).
$$
Off-diagonal elements in the sum in (\ref{posit}) will all be of the type
\begin{equation}
\pm \iint\limits_{A_j \times A_{j'} }\! f(x,y,u_{l,l'})\,dx \,dy 
\label{terms}
\end{equation}
and for each of the three combination $A_1\times A_1, A_1\times A_2$ and  $A_2\times A_2$
the number of i.i.d. terms of each type will be of order $n^2$, while the difference between 
the number of terms with opposite signs of each type will be at most $n/2$. Therefore, by the central 
limit theorem, the distribution of the left hand side of (\ref{posit}) converges weakly to some normal 
distribution and (\ref{posit}) can hold only if the variance of the terms in (\ref{terms}) is zero,
i.e. these terms are almost surely  constant. In particular,
$$
\iint\limits_{A_1\times A_{2} }\! f(x,y,u_{1,2})\,dx \,dy 
=
\iint\limits_{A_1 \times A_{2} }\! \bar{f}(x,y)\,dx \,dy
$$
with probability one. The same holds for some countable collection of sets $A_1\times A_2$
that generate the product $\sigma$-algebra on $[0,1]^2$ and this proves that for almost
all $z$ on $[0,1]$,  $f(x,y,z) = \bar{f}(x,y)$ for almost all $(x,y)$ on $[0,1]^2$. 
\qed\\
For simplicity of notations we will keep writing $f$ instead of $\bar{f}$ so that now
\begin{equation}
R_{l,l} = g(u_l)
\,\mbox{ and }\,
R_{l,l'} = f(u_l,u_{l'})
\,\mbox{ for }\,
l\not = l'
\label{fg3}
\end{equation}
is positive definite with probability one and $|f|, |g|\leq 1.$

\begin{lemma}\label{Lemphi}
If $R$ in (\ref{fg3}) is positive definite with probability one then
there exists a measurable map $\phi:[0,1]\to B$ where $B$ is the unit ball of $\ell^2$
such that
\begin{equation}
f(x,y) = \phi(x)\cdot \phi(y)
\label{fphi}
\end{equation}
almost surely on $[0,1]^2$.
\end{lemma}
\textbf{Remark.}
It is an important feature of the proof (similar to the argument in \cite{DS}) that the representation
(\ref{fphi}) of the off-diagonal elements is determined independently of the function $g$ that defines 
the diagonal elements. The diagonal elements play an auxiliary role in the proof of (\ref{fphi}) simply
through the fact that for some function $g$ the matrix $R$ in (\ref{fg3}) is positive definite. 
Once the representation (\ref{fphi}) is determined, the representation (\ref{repr}) will immediately follow.
\qed\\
\textbf{Proof.} Let us begin the proof with a simple observation that the fact that the matrix  (\ref{fg3})
is positive definite implies that $f(x,y)$ is a symmetric positive definite kernel on $[0,1]^2,$ 
\begin{equation}
\iint \! f(x,y) h(x) h(y) \,dx\,dy\geq 0 
\label{fint}
\end{equation}
for any $h\in L^2([0,1]).$ Since $(R_{l,l'})$ is positive definite, 
$n^{-2}\sum_{l,l'\leq n} R_{l,l'} h(u_l) h(u_{l'})\geq 0$ and since $|R_{l,l}|\leq 1,$ 
the diagonal terms $n^{-2}\sum_{l\leq n} R_{l,l} h(u_l)^2 \to 0$ a.s. as $n\to+\infty.$ 
Therefore, if we define
$$
S_n = \frac{2}{n(n-1)}\sum_{1\leq l<l'\leq n} f(u_l,u_{l'}) h(u_l) h(u_{l'}) 
$$
then $\liminf_{n\to+\infty} S_n \geq 0$ a.s. and (\ref{fint}) follows by the law of large numbers 
for $U$-statistics (Theorem 4.1.4 in \cite{decouple}), the proof of which we will recall for completeness. 
Namely, if we consider $\sigma$-algebra $\F_n = \sigma(u_{(1)},\ldots, u_{(n)}, (u_l)_{l>n})$ 
where $u_{(1)},\ldots, u_{(n)}$ are the order statistics of $u_1,\ldots, u_n$ then 
$(S_n,\F_n)$ is a reversed martingale and  $\bigcap_{n\geq 1} \F_n$ is trivial by 
the Hewitt-Savage zero-one law since it is in the tail $\sigma$-algebra of i.i.d. $(u_l)_{l\geq 1}$. 
Therefore, a.s.
$$
0\leq \lim_{n\to+\infty} S_n =  \e(S_2| \bigcap_{n\geq 1} \F_n) = \e S_2
$$
which proves (\ref{fint}). Since $f(x,y)$ is symmetric and in $L^2([0,1]^2),$ 
there exists an orthonormal sequence $(\varphi_l)$ in $L^2([0,1])$ such that 
(Theorem 4.2 in \cite{FA})
\begin{equation}
f(x,y) = \sum_{l \geq 1} \lambda_l\, \varphi_l(x) \varphi_l(y) 
\label{series}
\end{equation}
where the series converges in $L^2([0,1]^2).$ By (\ref{fint}), all $\lambda_l\geq 0$ and 
it is clear that now we would like to define $\phi$ in (\ref{fphi}) by 
\begin{equation}
\phi(x) = \bigl(\sqrt{\lambda_l}\,\varphi_l(x)\bigr)_{l\geq 1}: [0,1] \to \ell^2.
\label{phi}
\end{equation}
However, we still need to prove that the series in (\ref{series}) converges a.s. on $[0,1]^2$
and that $\sum_{l\geq 1}\lambda_l \,\varphi_l(x)^2\leq 1$ a.s. on $[0,1]$, i.e.  the map
$\phi$ is indeed into the unit ball of $\ell^2.$ For $m\geq 1,$ let $\DD_m$ be the $\sigma$-algebra 
on $[0,1]$ generated by the dyadic intervals $[i2^{-m},(i+1)2^{-m})$ so that 
$\DD_m\subset \DD_{m+1}$ and $\sigma\bigl(\bigcup_{m\geq 1} \DD_m\bigr)$ 
is the Borel $\sigma$-algebra on $[0,1].$ Let $\DD^m = \DD_m \otimes \DD_m$ 
be the product $\sigma$-algebra on $[0,1]^2$ and let
$$
f_n(x,y) = \sum_{1\leq l\leq n} \lambda_l\, \varphi_l(x) \varphi_l(y). 
$$
Since the series in (\ref{series}) converges in $L^2([0,1]^2),$  we can choose a subsequence 
$(n_j)$ such that the $L^2$-norm $\|f_{n_j} - f\|_2 \leq j^{-2}.$ Therefore,  
$\|\e(f_{n_j}| \DD^m) - \e(f| \DD^m)\|_2\leq j^{-2}$ for $m\geq 1$ and, by the Borel-Cantelli lemma, 
\begin{equation}
\lim_{j\to+\infty}\e(f_{n_j}| \DD^m) = \e(f| \DD^m) \,\,\rm a.s. 
\label{Dmart}
\end{equation}
Since $\DD^m$ is the product $\sigma$-algebra, $\e(f_{n_j}| \DD^m)$ is equal to
$$
\sum_{1\leq l\leq n_j} \lambda_l\, \e(\varphi_l | \DD_m)(x)\, \e(\varphi_l | \DD_m)(y)
$$
and for $x\in [i2^{-m},(i+1)2^{-m})$ we can define
$$
\e(\varphi_l | \DD_m)(x) = 2^m \!\int_{i2^{-m}}^{(i+1)2^{-m}}\!\!\! \varphi_l(s)\,ds.
$$ 
Therefore, for $x$ and $y$ in the same dyadic interval $[i2^{-m},(i+1)2^{-m})$,  
$$
\e(f_{n_j}| \DD^m)(x,y) = \sum_{1\leq l\leq n_j} \lambda_l\, (\e(\varphi_l | \DD_m))^2(x)
$$
and since $|\e(f|\DD^m)|\leq 1$, (\ref{Dmart}) implies that
$$
\lim_{j\to+\infty} \sum_{1\leq l\leq n_j} \lambda_l\, (\e(\varphi_l | \DD_m))^2(x) \leq 1 \,\, \rm a.s.
$$
The fact that all $\lambda_l\geq 0$ implies that for any $n\geq 1$
$$
\sum_{1\leq l\leq n} \lambda_l\, (\e(\varphi_l | \DD_m))^2(x) \leq 1\,\,\rm a.s.
$$  
By the martingale convergence theorem,  $\e(\varphi_l | \DD_m) \to \varphi_l$ a.s. as 
$m\to+\infty$ and, therefore, $\sum_{l\leq n} \lambda_l\, \varphi_l^2(x) \leq 1$ a.s.
Letting $n\to+\infty$  implies 
\begin{equation}
\phi(x)\cdot \phi(x) = \sum_{l \geq 1} \lambda_l\, \varphi_l^2(x) \leq  1 \,\,\rm a.s.
\label{phiphi}
\end{equation}
so the map $\phi$ in (\ref{phi}), indeed, maps $[0,1]$ into the unit ball of $\ell^2.$ Let us now show that
(\ref{fphi}) holds, i.e. the series in (\ref{series}) converges a.s. Given $n\geq 1,$ let us take 
$n_j\geq n$ and write $$ |f - f_n|\leq |f-f_{n_j}| + |f_{n_j}-f_n|.$$ The first term goes to zero a.s.
by the Borel-Cantelli lemma since $\|f_{n_j} - f\|_2 \leq j^{-2}$ and the second term can be bounded by
\begin{eqnarray*}
 |f_{n_j}(x,y)-f_n(x,y)| 
 &=& \Bigl|  \sum_{n<l\leq n_j} \! \lambda_l\, \varphi_l(x) \varphi_l(y)\Bigr|
 \\
 & \leq &
 \sum_{l>n} \! \lambda_l\, \varphi_l^2(x) \sum_{l>n} \! \lambda_l\, \varphi_l^2(y)
\end{eqnarray*}
and, by (\ref{phiphi}), also goes to zero a.s. as $n\to+\infty.$ Finally, the map $\phi$ is 
measurable since for any open ball $B_\eps(h)$ in $\ell^2$ of radius $\eps$ centered at $h,$
$\phi^{-1}(B_\eps(h))$ can be written as
$$
\sum_{l \geq 1} \lambda_l\, \varphi_l^2(x)
-2\sum_{l \geq 1} h_l \sqrt{\lambda_l}\, \varphi_l(x)
+h\cdot h< \eps
$$
and the left hand side is obviously a measurable function. This finishes the proof.
\qed\\
Lemma \ref{Lemphi} proves that if $(h_l,t_l)$ is an i.i.d. sequence from distribution 
$\eta = \lambda\circ (\phi,g)^{-1}$ on $\ell^2\times\Reals^+$ then the law of $R$ in (\ref{fg3}) 
coincides with the the law of 
\begin{equation}
\bigl( h_l\cdot h_{l'}(1-\delta_{l,l'}) + t_l\,\delta_{l,l'}\bigr)_{l,l'\geq 1}.
\label{reprbound}
\end{equation}
To prove (\ref{repr}) it remains to show that $\|h_l\|^2\leq t_l$ a.s. and define
$a_l = t_l - \|h_l\|^2.$ 
\begin{lemma}\label{Lemha}
The measure $\eta$  is concentrated on the set $\{(h,t) : \|h\|^2 \leq t\}.$
\end{lemma}
\textbf{Proof.} 
Suppose not. Then there exists $(h_0,t_0)$ such that $\|h_0\|^2 > t_0$ and such that the set
$$
A_\eps = B_\eps(h_0) \times (t_0-\eps,t_0+\eps)
$$ 
has positive measure $\eta(A_\eps)>0$ for all $\eps>0.$ 
Let us take $\eps>0$ small enough such that for any $(h_1,t_1), (h_2,t_2) \in A_\eps$ we have
$$
(h_1\cdot h_2)^2 > t_1\, t_2 + \eps.
$$ 
Since $\eta(A_\eps)>0,$ this contradicts the fact that for two independent copies 
$(h_1,t_1), (h_2,t_2)$ from distribution $\eta$ the matrix
$$
\left(
\begin{array}{cc}
t_1 & h_1\cdot h_2\\
h_1\cdot h_2 & t_2
\end{array}
\right)
$$
is positive definite with probability one.
\qed

\section{Unbounded case.}\label{SecUnbounded}

The idea of reducing the unbounded case to bounded one is briefly explained 
at the very end of the proof in \cite{DS} and here we will fill in the details.
Let us define a map $\Phi_N: M\to M$ such that for $\Gamma\in M,$
\begin{equation}
(\Psi_N(\Gamma))_{l,l'} = \Gamma_{l,l'}
\min\Bigl(\Bigl(\frac{N}{\Gamma_{l,l}}\Bigr)^{1/2},1\Bigr)
\min\Bigl(\Bigl(\frac{N}{\Gamma_{l',l'}}\Bigr)^{1/2},1\Bigr)
\label{PsiN}
\end{equation}
and, in particular, $(\Psi_N(\Gamma))_{l,l}=\min(N,\Gamma_{l,l}).$ Define a map
$\psi_N: \ell^2\times \Reals^+ \to \ell^2\times \Reals^+$ by
\begin{equation}
\psi_N(h,t)=\Bigl( h\min\Bigl(\Bigl(\frac{N}{t}\Bigr)^{1/2},1\Bigr),
\min(N,t) \Bigr).
\label{psiN}
\end{equation}
Let us make two simple observations that follow from the definitions (\ref{PsiN}) and (\ref{psiN}):
\begin{enumerate}
\item[(a)]{ if $\Gamma$ is a Gram-de Finetti matrix then $\Psi_N(\Gamma)$ is also a Gram-de Finetti
matrix uniformly bounded by $N$ with probability one;
}
\item[(b)]{ if $\Gamma$ is a Gram-de Finetti matrix generated as in (\ref{reprbound}) by an i.i.d. 
sample $(h_l, t_l)$ from distribution $\nu$ on $\ell^2\times \Reals^+$ then $\Psi_N(\Gamma)$ 
is generated by an i.i.d. sample $(\psi_N(h_l,t_l))$ from distribution $\nu\circ \psi_N^{-1}.$
}
\end{enumerate}
Consider a Gram-de Finetti matrix $R.$ Since $\Psi_N(R)$ is uniformly bounded, the results of
Section \ref{SecBounded} imply that it can be generated as in (\ref{reprbound}) by an i.i.d.
sample from some measure $\eta_N$ on $\ell^2\times\Reals^+.$ Since 
$\lim_{N\to+\infty}\Psi_N(R)=R$ a.s., this indicates that $R$ should be generated by an i.i.d. 
sample from distribution $\eta$ defined as a limit of $\eta_N.$ However, to ensure that 
this limit exists we first need to redefine the sequence $(\eta_N)$ in a consistent way.
For this, we will need to use the fact that a measure $\eta$ in the representation (\ref{reprbound}) 
is unique up to an orthogonal transformation of its marginal on $\ell^2.$

\begin{lemma}\label{Lemq}
If $(h_l,t_l)$ and $(h_l',t_l')$ are i.i.d. samples from distributions $\eta$ and $\eta'$
on $\ell^2\times\Reals^+$ correspondingly and 
\begin{equation}
\bigl( h_l\cdot h_{l'}(1-\delta_{l,l'}) + t_l\,\delta_{l,l'}\bigr)
\stackrel{d}{=}
\bigl( h_l'\cdot h_{l'}'(1-\delta_{l,l'}) + t_l'\,\delta_{l,l'}\bigr)
\label{repreq}
\end{equation}
then there exists a unitary operator $q$ on $\ell^2$ such that $\eta  = \eta' \circ (q, \id)^{-1}.$
\end{lemma}
\textbf{Proof.}
Let us begin by showing that the values $h_l\cdot h_{l}$ and $h_l'\cdot h_{l}'$ can be reconstructed 
almost surely from the matrices (\ref{repreq}). Consider a sequence $(g_l)$ on $\ell^2$ such that 
$\|g_l\|^2 = t_l$ and $g_l\cdot g_{l'} = h_l\cdot h_{l'}$ for all $l<l'$.  Without loss of generality, 
let us assume that
$$
g_l = h_l + \sqrt{t_l-\|h_l\|^2} e_l
$$
where $(e_l)$ is an orthonormal sequence orthogonal to the closed span of $(h_l)$
(if necessary, we identify $\ell^2$ with $\ell^2\oplus \ell^2$ to choose the sequence $(e_l)$).
Since $(h_l)$ is an i.i.d. sequence from the marginal $\mu$ of measure $\eta$ on $\ell^2$, 
with probability one there are elements in the sequence $(h_l)_{l\geq 2}$ arbitrarily close to $h_1$ 
and, therefore, the length of the orthogonal projection of $h_1$ onto the closed span of 
$(h_l)_{l\geq 2}$ is equal to $\|h_1\|.$ As a result, the length of the orthogonal projection of $g_1$ 
onto the closed span of $(g_l)_{l\geq 2}$ is also equal to $\|h_1\|$ which means that we reconstructed 
$\|h_1\|$ from the first matrix in (\ref{repreq}). Therefore, (\ref{repreq}) implies that
\begin{equation}
\bigl( (h_l\cdot h_{l'}), (t_l) \bigr)
\stackrel{d}{=}
\bigl((h_l'\cdot h_{l'}'),(t_l')\bigr).
\label{hsfull}
\end{equation}
Given $(h_l\cdot h_{l'})$ and $(h_l'\cdot h_{l'}')$, we can now construct sequences $(x_l)$ and $(x_l')$
isometric to $(h_l)$ and $(h_l')$ in some pre-determined way, for example, by choosing $x_l$ and 
$x_l'$ to be in the span of the first $l$ elements of some fixed orthonormal basis. Then there exist
(random) unitary operators $U = U((h_l)_{l\geq 1})$ and $U'=U'((h_l')_{l\geq 1})$ on $\ell^2$ such that 
\begin{equation}
x_l = U h_l \,\,\mbox{ and }\,\, x_l' = U' h_{l}'.
\label{UV}
\end{equation} 
By the strong law of large number for empirical measures (Theorem 11.4.1 in \cite{Dudley})
$$
\frac{1}{n}\sum_{1\leq l\leq n} \delta_{(h_l,t_l)} \to \eta
\,\,\mbox{ and }\,\,
\frac{1}{n}\sum_{1\leq l\leq n}\delta_{(h_l',t_l')} \to \eta' 
$$
weakly almost surely, and therefore, (\ref{UV}) implies that
$$
\frac{1}{n}\sum_{1\leq l\leq n}\delta_{(x_l,t_l)} \to \eta \circ (U,\id)^{-1} 
\,\,\mbox{ and }\,\,
\frac{1}{n}\sum_{1\leq l\leq n}\delta_{(x_l',t_l')} \to \eta' \circ (U',\id)^{-1} 
$$
weakly almost surely. Therefore, since $(x_l,t_l)$ and $(x_l',t_l')$ have the same distribution 
by (\ref{hsfull}), $\eta \circ (U,\id)^{-1}$ and $\eta' \circ (U',\id)^{-1}$ have the same distribution
on the space of all probability distributions on $\ell^2\times \Reals^+$ with the topology of weak
convergence. This implies that there exist non-random unitary operators $U$ and $U'$ 
such that $\eta \circ (U,\id)^{-1} = \eta' \circ (U',\id)^{-1}$ and taking $q=U^{-1}U'$ finishes the proof.
\qed\\
Using Lemma \ref{Lemq}, we will now construct a "consistent" sequence of laws $(\eta_N)$ 
recursively as follows. Suppose that the measure $\eta_N$ that generates $\Psi_N(R)$ 
as in (\ref{reprbound}) has already been defined. Suppose now that $\Psi_{N+1}(R)$ is generated 
by an i.i.d. sample $(h_l,t_{l})$ from some measure $\eta_{N+1}$. Since
\begin{equation}
\Psi_{N}(\Psi_{N+1}(R)) = \Psi_{N}(R),
\label{Psis}
\end{equation}
observation (b) above implies that  $\Psi_{N}(R)$ can also be generated by $\psi_{N}(h_l,t_{l})$
from measure $\eta_{N+1}\circ \psi_N^{-1}.$ Lemma \ref{Lemq} then implies that there exists
a unitary operator $q$ on $\ell^2$ such that
$$
\eta_N = (\eta_{N+1}\circ \psi_N^{-1}) \circ (q,\id)^{-1} =
(\eta_{N+1} \circ (q,\id)^{-1})\circ \psi_N^{-1}
$$
since $\psi_N$ and $(q,\id)$ obviously commute. We now redefine $\eta_{N+1}$ to be equal to 
$\eta_{N+1} \circ (q,\id)^{-1}$. Clearly, $\Psi_{N+1}(R)$ is still generated by an i.i.d. sequence 
from this new measure $\eta_{N+1}$ and in addition we have 
\begin{equation}
\eta_N = \eta_{N+1} \circ \psi_N^{-1}.
\label{consist}
\end{equation}
Let $A_N:=\ell^2\times [0,N).$ Since $\psi_N(h,t)\in A_N$ if and only if $(h,t)\in A_N$ and 
$\psi_N(h,t) = (h,t)$ on $A_N$, the consistency condition (\ref{consist}) implies that the restrictions
of measures $\eta_N$ and $\eta_{N+1}$ to $A_N$ are equal. Therefore, $\eta_N$ converges
in total variation to
$ \eta = \sum \eta_N\!\! \downharpoonright_{A_N\setminus A_{N-1}}$. 
Since $\psi_N = \psi_N\circ \psi_{N'}$ for $N\leq N'$, (\ref{consist}) implies that
$\eta_N = \eta_{N'} \circ \psi_N^{-1}$ and since $\psi_{N}$ is continuous, letting $N'\to+\infty$
gives $\eta_N = \eta \circ \psi_N^{-1}.$ Finally, letting $N\to+\infty$ proves that $R$ is generated 
by  an i.i.d. sample  from $\eta$ which proves representation (\ref{reprbound})
in the unbounded case, and Lemma \ref{Lemha} again implies (\ref{repr}).
\qed

\end{document}

%% file: invlat.tex
\font\tencmmib=cmmib10 \skewchar\tencmmib '60
\newfam\cmmibfam
\textfont\cmmibfam=\tencmmib

\def\bbox{\quad\hbox{\vrule \vbox{\hrule \vskip2pt \hbox{\hskip2pt
\vbox{\hsize=1pt}\hskip2pt} \vskip2pt\hrule}\vrule}}
\def\lessim{\ \lower4pt\hbox{$
\buildrel{\displaystyle <}\over\sim$}\ }
\def\gessim{\ \lower4pt\hbox{$\buildrel{\displaystyle >}
\over\sim$}\ }

\def\DD{{\cal D} }

\def\eps{\varepsilon}

\def\go0{\to 0}

\def\leftitem#1{\item{\hbox to\parindent{\enspace#1\hfill}}}

\def\qed{\hfill\break\rightline{$\bbox$}}

\def\sg{\sigma}

\def\sg2{\sigma^2}

\def\__{_{\infty}}